\documentclass[11pt]{amsart}
\usepackage[english]{babel}
\usepackage{graphicx}
\newtheorem{theorem}{Theorem}[section]

\newtheorem{lemma}[theorem]{Lemma}

\theoremstyle{definition}

\theoremstyle{remark}
\newtheorem{remark}{Remark}
\numberwithin{equation}{section}
\newcommand{\real}{\mathbb R}
\def\natu{\mathbb N}
\def\dis{\displaystyle}
\def\sobol{H^{1}(0,L)}

\def\landan{\Lambda_n}
\def\landaene{\lambda_n}

\begin{document}

\title[Lyapunov inequalities at higher eigenvalues]
{Lyapunov inequalities for Neumann boundary conditions at higher eigenvalues}%
\author{A. Ca\~{n}ada}
\author{S. Villegas}
\thanks{The authors have been supported by the Ministry of Education and
Science of Spain (MTM2005-01331)}
\address{Departamento de An\'{a}lisis
Matem\'{a}tico, Universidad de Granada, 18071 Granada, Spain.}
\email{acanada@ugr.es, svillega@ugr.es}


\begin{abstract}
This paper is devoted to the study of Lyapunov-type inequality for
Neumann boundary conditions at higher eigenvalues. Our main result
is derived from a detailed analysis about the number and
distribution of zeros of nontrivial solutions and their first
derivatives, together with the use of suitable minimization
problems. This method of proof allows to obtain new information on
Lyapunov constants. For instance, we prove that as in the
classical result by Lyapunov, the best constant is not attained.
Additionally, we exploit the relation between Neumann boundary
conditions and disfocality to provide new nonresonance conditions
at higher eigenvalues.
\end{abstract}
\maketitle
\section{Introduction}

\noindent The classical $L^1$ Lyapunov inequality for Neumann
boundary problem
\begin{equation}\label{n1}
u''(x) + a(x)u(x) = 0, \ x \in (0,L), \ u'(0) = u'(L) = 0
\end{equation}
states that if
\begin{equation}\label{2806072}
a \in L^1(0,L)\setminus \{0\},\  \displaystyle \int_{0}^{L} a(x) \
dx \geq 0 \end{equation} is such that (\ref{n1}) has nontrivial
solutions, then $\displaystyle \int_{0}^{L} a^{+}(x) \ dx
> 4/L$, where $a^{+}(x) = \max \{ a(x),0 \}$ (\cite{ha},
\cite{huyo}). In \cite{camovi} and \cite{zhang} the authors
generalize this result providing, for each $p$ with $1 \leq p \leq
\infty,$ optimal necessary conditions for  boundary value problem
(\ref{n1}) to have nontrivial solutions, given in terms of the
$L^p$ norm of the function $a^+.$ In particular, if $p = \infty,$
it is proved that (\ref{n1}) has only the trivial solution if
function $a$ satisfies
\begin{equation}\label{n1b}
 a\in L^\infty (0,L) \setminus \{ 0 \},\
\displaystyle\int_0^L a \geq 0, \ \ a^+ \prec \pi^2/L^2,
\end{equation}
where for $c,d \in L^1 (0,L),$ we write $c \prec d$ if $c(x) \leq
d(x)$ for a.e. $x \in [0,L]$ and $c(x) < d(x)$ on a set of
positive measure.  This is a very well known result which is
usually called the nonuniform nonresonance condition with respect
to the two first eigenvalues $\lambda_0 = 0 $ and $\lambda_1 =
\pi^2/L^2$ of the eigenvalue problem
\begin{equation}\label{n2}
 u''(x) + \lambda u(x) = 0, \ x \in (0,L), \ u'(0) = u'(L) =
0 \end{equation} (see \cite{mawa}, \cite{mawa2} and
\cite{mawawi}). From this point of view, it may be affirmed that
the nonuniform nonresonance condition (\ref{n1b}) is in fact the
$L_\infty$ Lyapunov inequality at the two first eigenvalues
$\lambda_0$ and $\lambda_1$.

On the other hand, the set of eigenvalues of (\ref{n2}) is given
by $\lambda_n = n^2 \pi^2/L^2, \ n \in \natu \cup \{ 0 \}$ and by
using a general result due to Dolph \cite{dolp}, it can be proved
that, if for some $n \geq 1$ function $a$ satisfies
\begin{equation}\label{n1bb}
\lambda_n \prec a \prec \lambda_{n+1}
\end{equation}
then (\ref{n1}) has only the trivial solution (see \cite{mawa3},
Lemma 2.1, for some generalizations of (\ref{n1bb}) to more
general boundary value problems). It is clear that condition
(\ref{n1bb}) can not be obtained from $L_p$ Lyapunov inequalities
given in \cite{camovi} and \cite{zhang}.

Previous observations motivate this article where, for any given
natural number $n \geq1 $ and function $a$ satisfying $\lambda_n
\prec a, $we obtain the $L_1$ Lyapunov inequality (the case of
$L_p$ with $1 <p<\infty$ presents some special particularities and
will be considered in a forthcoming paper). In particular we
prove, as in the classical Lyapunov inequality, that the best
constant is not attained for any value of $n.$ To the best of our
knowledge this result is new if $n \geq 1$. In the $L^\infty$
case, Lyapunov inequality is exactly (\ref{n1bb}) and in this
sense, it is natural to say that this paper deals with Lyapunov
inequalities at higher eigenvalues.

One of the main results of our paper is given by Lemma
\ref{280607l1} below where we discuss in detail the number and
distribution of zeros of $u$ and $u',$ where $u$ is any nontrivial
solution of the linear boundary value problem (\ref{n1}).

In the second section we study the $L^1$ Lyapunov inequality when
$\landaene \prec a$. The case where function $a$ satisfies the
condition $ A \leq a(x) \leq B, \ \mbox{a.e. in} \ (0,L) $ where
$\lambda_k < A < \lambda_{k+1} \leq B$ for some $k \in \natu \cup
\{ 0 \}, $ has been considered in \cite{yong}. In this paper the
authors use Optimal Control theory methods, specially Pontryagin's
maximum principle.

In the last section we use the natural relation between Neumann
boundary conditions and disfocality, given by Lemma
\ref{280607l1}, to obtain new results on the existence and
uniqueness of solutions for linear resonant problems with Neumann
boundary conditions. We use $L^1$ and $L^\infty$ Lyapunov
constants. For example, by using Lemma \ref{280607l1} and the
$L^\infty$ Lyapunov inequality, we can prove (see Theorem
\ref{t2211} below) that if
\begin{equation}\label{211107}
\begin{array}{c}
a \in L^\infty (0,L), \ \landaene \prec a \
\mbox{and} \ \exists \ \ 0 = y_0 < y_1 < \ldots < y_{2n+1} < y_{2n+2} = L: \\ \\
\max\ _{0 \leq i \leq 2n+1}  \{ (y_{i+1}-y_i)^2 \Vert a
\Vert_{L^\infty (y_i,y_{i+1})} \} \leq \pi^2/4
\end{array}
\end{equation}
 and, in addition, $a$ is not the
constant $\pi^2 /4(y_{i+1}-y_i)^2,$ at least in one of the
intervals $[y_i,y_{i+1}], \ 0 \leq i \leq 2n+1,$ then we obtain
that (\ref{n1}) has only the trivial solution (this kind of
functions $a$ are usually named $2(n+1)$-step potentials).

Previous hypothesis is optimal in the sense that if $a$ is the
constant $\pi^2 /4(y_{i+1}-y_i)^2$ in each one of the intervals
$(y_i,y_{i+1}), \ 0 \leq i \leq 2n+1,$ then (\ref{n1}) has
nontrivial solutions (see Remark \ref{r221107} in section 3). If
$y_i = \frac{iL}{2(n+1)},\ 0 \leq i \leq 2n+2,$ we have the so
called non-uniform non-resonance conditions at higher eigenvalues
(\cite{dolp}, \cite{mawa3}) but if for instance, $y_{j+1}-y_{j} <
\frac{L}{2(n+1)},$ for some $j, \ 0 \leq j \leq 2n+1,$ function
$a$ can satisfies $\Vert a \Vert_{L^\infty (y_j,y_{j+1})} =
\frac{\pi^2}{4(y_{j+1}-y_j)^2}$ (which is a quantity greater than
$\lambda_{n+1} = \frac{(n+1)^2\pi^2}{L^2}$) as long as $a$
satisfies (\ref{211107}) for each $i \neq j.$

Additionally, such as it has been done in \cite{camovi},
\cite{camovi2}, \cite{cavi} and \cite{yong}, the linear study can
be combined with Schauder fixed point theorem to provide new
conditions about the existence and uniqueness of solutions for
resonant nonlinear problems (see Theorem \ref{t2bis} below). Also,
we may deal with other boundary value problems. Finally, one can
expect that some results hold true in the case of Neumann boundary
value problem for partial differential equations
\begin{equation}\label{pde} \Delta u(x) + a(x)u(x) = 0, \ x \in \Omega, \ \dis
\frac{\partial u(x)}{\partial n} = 0, \ x \in \partial \Omega
\end{equation} where $\Omega$ is a bounded and regular domain in
$\real^{N},$ but here the role played by the dimension $N$ may be
important (see \cite{camovi2}).

\section{Lyapunov inequality at higher eigenvalues}

If $n \in \natu$ is fixed, we introduce the set $\landan$ as

\begin{equation}\label{2806073}
\landan = \{ a \in L^{1}(0,L): \lambda_n \prec a \ \mbox{and} \
(\ref{n1}) \ \mbox{has nontrivial solutions} \ \}
\end{equation}

Here $u \in \sobol,$ the usual Sobolev space. If we define the
number
\begin{equation}\label{1009071}
\beta_{1,n} \equiv \inf_{a \in \landan } \ \Vert a - \landaene
\Vert_{L^1 (0,L)}
\end{equation} the main result of this section is the following.

\begin{theorem}\label{t1}
$$
\beta_{1,n}  = \frac{2\pi n(n+1)}{L} \cot \frac{\pi n}{2(n+1)}
$$
Moreover $\beta_{1,n}$ is not attained.
\end{theorem}

\begin{proof}
It is based on some lemmas. In the first one we do a careful
analysis about the number and distribution of zeros of the
nontrivial solutions $u$ of (\ref{n1}). Since $a \in \landan,$ it
is clear that between two consecutive zeros of the function $u$
there must exists a zero of the function $u'$ and between two
consecutive zeros of the function $u'$ there must exists a zero of
the function $u.$ More precisely, we have the following lemma.

\begin{lemma}\label{280607l1}
Let $a \in \landan$ be given and $u$ any nontrivial solution of
(\ref{n1}). If the zeros of $u'$ in $[0,L]$ are denoted by $0 =
x_0 < x_2 < \ldots < x_{2m} = L$ and the zeros of $u$ in $(0,L)$
are denoted by $ x_1 < x_3 < \ldots < x_{2m-1}, $ then:
\begin{enumerate}
\item $x_{i+1} - x_{i} \leq \frac{L}{2n}, \ \forall \ i: \ 0 \leq
i \leq 2m-1.$ Moreover, at least one of these inequalities is
strict. \item $m \geq n+1.$ Moreover, any value $m\geq n+1$ is
possible. \item Let $i, \ 0 \leq i \leq 2m-1,$ be given. Then,
functions $a$ and $u$ satisfy
\begin{equation}\label{02071}
\Vert a - \landaene \Vert_{L^1 (x_i,x_{i+1})} \geq
\frac{\int_{x_i}^{x_{i+1}} u'^2 - \landaene \int_{x_i}^{x_{i+1}}
u^2 }{u^2(x_{i+1})}, \ \mbox{if} \ i \ \mbox{is odd}
\end{equation}
and
\begin{equation}\label{02072}
\Vert a - \landaene \Vert_{L^1 (x_i,x_{i+1})} \geq
\frac{\int_{x_i}^{x_{i+1}} u'^2 - \landaene \int_{x_i}^{x_{i+1}}
u^2 }{u^2(x_i)}, \ \mbox{if} \  i \ \mbox{is even}
\end{equation}
\end{enumerate}
\end{lemma}
\begin{proof}
Let $i, \ 0 \leq i \leq 2m-1,$ be given. Then, function $u$
satisfies either the problem \begin{equation}\label{2806076}
u''(x) + a(x)u(x) = 0, \ x \in (x_i,x_{i+1}), \ \ u(x_i) = 0, \
u'(x_{i+1}) = 0, \end{equation} or the problem
\begin{equation}\label{2806077} u''(x) + a(x)u(x) = 0, \ x \in (x_i,x_{i+1}),
\ \ u'(x_i) = 0, \ u(x_{i+1}) = 0. \end{equation} Let us assume
the first case. The reasoning in the second case is similar. Note
that $u$ may be chosen such that $u(x)
> 0, \ \forall \ x  \in  (x_i,x_{i+1}).$ Let us denote by $\mu_1 ^i$ and
$\varphi_1^i,$ respectively, the principal eigenvalue and
eigenfunction of the eigenvalue problem
\begin{equation}\label{2806074} v''(x) + \mu v(x) = 0, \ x \in
(x_i,x_{i+1}), \ v(x_i) = 0, \ v'(x_{i+1}) = 0. \end{equation} It
is known that
\begin{equation}\label{2806075}
\mu_1 ^i = \frac{\pi^2}{4(x_{i+1}-x_i)^2}, \ \varphi_1 ^i (x) =
\sin \frac{\pi (x-x_i)}{2(x_{i+1}-x_i)} \end{equation} Choosing
$\varphi_1 ^i$ as test function in the weak formulation of
(\ref{2806076}) and $u$ as test function in the weak formulation
of (\ref{2806074}) for $\mu = \mu_1 ^i$ and $v =\varphi_1^i,$ we
obtain
\begin{equation}\label{2806078}
\int_{x_i}^{x_{i+1}} (a(x) - \mu_1^i)u\varphi_1^i (x) \ dx = 0.
\end{equation}
Then, if $x_{i+1} - x_{i} > \frac{L}{2n},$ we have
$$
\mu_1 ^i = \frac{\pi^2 L^2}{4(x_{i+1}-x_i)^2L^2} < \frac{n^2
\pi^2}{L^2} = \landaene \leq a(x), \ \mbox{a.e. in} \
(x_i,x_{i+1})
$$
which is a contradiction with (\ref{2806078}). Consequently,
$x_{i+1} - x_i \leq \frac{L}{2n}, \forall \ i: \ 0 \leq i \leq
2m-1.$ Also, since $\lambda_n \prec a$ in the interval $(0,L)$, we
must have $\lambda_n \prec a$ in some subinterval $(x_j,x_{j+1})$.
If $x_{j+1} - x_{j}= \frac{L}{2n},$ it follows $\mu_1 ^j \prec a$
in $(x_j,x_{j+1})$ and this is again a contradiction with
(\ref{2806078}). These reasonings complete the first part of the
lemma. For the second one, let us observe that
$$
L = \sum_{i=0}^{2m-1} (x_{i+1}-x_i) < 2m \frac{L}{2n}
$$
In consequence, $m > n.$ Also, note that for any given natural
number $q \geq n+1,$ function $a(x) \equiv \lambda_{q}$ belongs to
$\landan$ and for function $u(x) = \cos \frac{q\pi x}{L},$ we have
$m = q.$

Lastly, if $i,$ with $ \ 0 \leq i \leq 2m-1$ is given and $u$
satisfies (\ref{2806076}), then
$$
\begin{array}{c}
\int_{x_i}^{x_{i+1}} u'^2(x) = \int_{x_i}^{x_{i+1}} a(x)u^2(x) =
\\ \\
\int_{x_i}^{x_{i+1}} (a(x)-\landaene) u^2 (x) +
\int_{x_i}^{x_{i+1}} \landaene u^2 (x)
\end{array}
$$
Therefore, $$ \int_{x_i}^{x_{i+1}} u'^2(x) - \landaene
\int_{x_i}^{x_{i+1}} u^2(x) \leq \Vert a - \landaene \Vert_{L^1
(x_i,x_{i+1})} \Vert u^2 \Vert_{L^\infty (x_i,x_{i+1})}
$$
Since $u'$ has no zeros in the interval $(x_i,x_{i+1})$ and
$u(x_i) = 0,$ we have $\Vert u^2 \Vert_{L^\infty(x_i,x_{i+1})} =
u^2 (x_{i+1}).$ This proves the third part of the lemma when $u$
satisfies (\ref{2806076}). The reasoning is similar if $u$
satisfies (\ref{2806077}).
\end{proof}

\begin{lemma}\label{l2}
Assume that $a < b$ and $0 < M \leq \frac{\pi^2}{4(b-a)^2}$ are
given real numbers. Let $H = \{ u \in H^1 (a,b) : u(a) = 0, u(b)
\neq 0 \}.$ If $J: H \rightarrow \real$ is defined by
\begin{equation}\label{02073}
J(u) = \frac{\int_{a}^{b} u'^2 - M \int_{a}^{b} u^2 }{u^2(b)}
\end{equation}
and $m \equiv \inf_{u \in H} \ J(u),$ then $m$ is attained.
Moreover
\begin{equation}
m = M^{1/2}\cot (M^{1/2}(b-a))
\end{equation}
and if $u \in H,$ then $J(u) = m \Longleftrightarrow u(x) = k
\frac{\sin (M^{1/2}(x-a))}{\sin (M^{1/2}(b-a))}$ for some non zero
constant $k.$
\end{lemma}
\begin{proof}
Remember that $\delta_1 = \frac{\pi^2}{4(b-a)^2}$ is the principal
eigenvalue of the eigenvalue problem $v''(x) + \delta v(x) = 0, \
v(a) = 0, \ v'(b) = 0$ with associated eigenfunction $w(x) = \sin
\frac{\pi(x-a)}{2(b-a)}.$ Therefore, if $M
=\frac{\pi^2}{4(b-a)^2},$ $m = 0$ and it is attained at function
$w.$

If $M < \delta_1 = \frac{\pi^2}{4(b-a)^2},$ there exists some
positive constant $c$ such that
\begin{equation}\label{02074}
\int_a^b u'^2 - M \int_a^b u^2 \geq c \int_a^b u'^2, \ \forall \ u
\in H.
\end{equation}
If  $\{ u_n \} \subset H$ is a minimizing sequence for $J$, since
the sequence $\{ k_n u_n \}, \ k_n \neq 0,$ is also a minimizing
sequence, we can assume without loos of generality that $u_n (b)
=1.$ From (\ref{02074}) we deduce that $\int_a^b u_n '^2 $ is
bounded. So, we can suppose, up to a subsequence, that $u_n
\rightharpoonup u_0$ in $H^1 (a,b)$ and $u_n \rightarrow u_0$ in
$C[a,b]$ (with the uniform norm). The strong convergence in
$C[a,b]$ gives us $u_0 (b) = 1.$ The weak convergence in $H$
implies $J(u_0) \leq \liminf J(u_n) = m.$ Then $u_0$ is a
minimizer.

Since $J(u_0) = \min \{J(v): v \in H^1 (a,b), v(a) = 0, v(b) =1
\},$ Lagrange multiplier Theorem implies that there are real
numbers $\alpha_1, \alpha_2$ such that
$$
2 \int_a^b u_0'v' - 2M \int_a^b u_0 v - \alpha_1 v(b) - \alpha_2
v(a) = 0, \ \forall v \in H^1(a,b).
$$
In particular,
$$
 \int_a^b u_0 'v' - M \int_a^b u_0 v =
0, \ \forall v \in H^1(a,b): v(a) = v(b) = 0.
$$
We conclude that $u_0$ satisfies the problem
\begin{equation}\label{02075}
u_0 ''(x) + Mu_0 (x) = 0, \ x \in (a,b), \ u_0 (a) = 0, \ u_0 (b)
= 1.
\end{equation}
Note that since $M < \frac{\pi^2}{(b-a)^2},$ (\ref{02075}) has a
unique solution, which is given by
\begin{equation}
u_0 (x) = \frac{\sin (M^{1/2}(x-a))}{\sin (M^{1/2}(b-a))}.
\end{equation}
Finally, an elementary calculation gives $J(u_0) = M^{1/2}\cot
(M^{1/2}(b-a)).$ This proves the lemma.
\end{proof}
Now, we combine Lemma \ref{280607l1} and Lemma \ref{l2} to obtain
the following result.

\begin{lemma}\label{l3}
Let $a \in \landan$ be given and $u$ any nontrivial solution of
(\ref{n1}). If the zeros of $u'$ are denoted by $0 = x_0 < x_2 <
\ldots < x_{2m} = L$ and the zeros of $u$ are denoted by $ x_1 <
x_3 < \ldots < x_{2m-1}, $ then:
\begin{equation}
\Vert a - \landaene \Vert_{L^1 (0,L)} \geq \frac{n\pi}{L}
\sum_{i=0}^{2m-1} \cot (\frac{n\pi}{L} (x_{i+1} -x_i))
\end{equation}
\end{lemma}

Previous reasoning motivates the study of a special minimization
problem given in the following lemma.

\begin{lemma}\label{l4}
Given any $r\in\natu$ and $S\in\real^+$ satisfying $r \pi > 2S$,
let
$$
Z = \{ z= (z_0,z_1,\ldots,z_{r-1}) \in (0,\pi/2]^r: \
\sum_{i=0}^{r-1} z_i = S \}
$$
If $F: Z \rightarrow \real$ is defined by
$$
F(z) = \sum_{i=0}^{r-1} \ \cot \ z_i,
$$
then $ \displaystyle \inf_{z \in Z} \ F(z)$ is attained and its
value is $ r \cot \frac{S}{r}.$ Moreover, $z \in Z$ is a minimizer
if and only if $z_i=\frac{S}{r}, \ \forall \ 0 \leq i \leq r-1.$
\end{lemma}
\begin{proof}
Let us observe that $\forall \ z \in Z,\  \cot \ z_i \geq 0, \ 0
\leq i \leq r-1.$ Moreover, if $z_i \rightarrow 0^+$ for some $0
\leq i \leq r-1$, then $\cot \ z_i \rightarrow + \infty.$ Also,
since $r \pi > 2 S,$ if $z \in Z$ is such that $z_i = \pi/2,$ for
some $0 \leq i \leq r-1,$ then there must exist some $0 \leq j
\leq r-1$ such that $z_j < \pi/2.$ Let us choose the point $z^*
\in Z$ defined (for $\delta
>0$ sufficiently small) as $z^*_k = z_k,$ if $k \neq i$ and $k
\neq j,$ $z^*_i = \displaystyle \frac{\pi}{2} - \delta, \ z^*_j =
z_j + \delta.$ An elementary calculation shows
$$
F(z^*) - F(z) = \displaystyle \frac{\cot z_j (1-\cot z_j \cot
\delta)}{\cot \delta (\cot z_j + \cot \delta)}
$$
which is a negative number for $\delta$ sufficiently small.
Consequently, there exits a sufficiently small positive constant
$\varepsilon_1$ such that
$$ \inf_{z \in Z} \ F(z) = \min_{z \in
[\varepsilon_1,\frac{\pi}{2}]^{r}} \ F(z) = \min_{z \in
(\varepsilon_1,\frac{\pi}{2} )^{r}} \ F(z)$$ Then, if $z \in Z$ is
any minimizer of $F,$ Lagrange multiplier Theorem implies that
there is $\lambda \in \real$ such that
$$
\frac{-1}{\sin^2 z_i} + \lambda = 0,\ 0 \leq i \leq r-1, \ \
\sum_{i=0}^{r-1} z_i = S.
$$
We conclude $z_i = \frac{S}{r}, \ 0 \leq i \leq r-1$ and the lemma
is proved.
\end{proof}

From two previous lemmas, we obtain the following one.

\begin{lemma}\label{l5}
\begin{equation}\label{03071}
\beta_{1,n}  \ \geq \frac{n\pi}{L} 2(n+1)\cot \frac{n\pi}{2(n+1)}.
\end{equation}
\end{lemma}
\begin{proof}
Let $a \in \landan$ be given and $u$ any nontrivial solution of
(\ref{n1}). If the zeros of $u'$ are denoted by $0 = x_0 < x_2 <
\ldots < x_{2m} = L$ and the zeros of $u$ are denoted by $ x_1 <
x_3 < \ldots < x_{2m-1}, $ then we obtain from Lemma \ref{l3} and
Lemma \ref{l4} (with $r=2m$, $S=n\pi$ and
$z_i=\frac{n\pi}{L}(x_{i+1}-x_i)$)
\begin{equation}
\Vert a - \landaene \Vert_{L^1 (0,L)} \geq \frac{n\pi}{L}
\sum_{i=0}^{2m-1} \cot (\frac{n\pi}{L} (x_{i+1} -x_i)) \geq
\frac{n\pi}{L} 2m \cot \frac{n\pi}{2m}.
\end{equation}
Finally, taking into account the property
$$
\mbox{The function}\ 2m \cot \frac{n\pi}{2m}\ \mbox{is strictly
increasing with respect to} \ m \eqno(P)$$ and that $m \geq n+1,$
we deduce (\ref{03071}).
\end{proof}

In the next lemma, we define a minimizing sequence for
$\beta_{1,n}.$

\begin{lemma}\label{l7}
Let $ \varepsilon > 0$ be sufficiently small. Let us define the
function $u_\varepsilon : [0,L] \rightarrow \real$ by
\begin{equation}
u_\varepsilon (x) = \left \{
\begin{array}{l}
-\sin (\frac{n\pi}{L}(x-\frac{L}{2(n+1)})) +
\frac{n\pi}{L}\frac{(x-\varepsilon)^3}{3\varepsilon^2}\cos
(\frac{n\pi}{2(n+1)}), \ \ \mbox{if}  \ \ 0 \leq x \leq \varepsilon, \\
\\
-\sin (\frac{n\pi}{L}(x-\frac{L}{2(n+1)})), \ \ \mbox{if} \ \
\varepsilon \leq x \leq \frac{L}{2(n+1)}, \\ \\
-u_\varepsilon (\frac{2L}{2(n+1)} -x), \ \ \mbox{if} \ \
\frac{L}{2(n+1)} \leq x \leq \frac{2L}{2(n+1)}, \\ \\
u_\varepsilon (\frac{4L}{2(n+1)}-x), \ \ \mbox{if} \ \
\frac{2L}{2(n+1)} \leq x \leq \frac{4L}{2(n+1)}, \\ \\
-u_\varepsilon (\frac{6L}{2(n+1)}-x), \ \ \mbox{if} \ \
\frac{4L}{2(n+1)} \leq x \leq \frac{6L}{2(n+1)},  \\ \\
\ldots
\end{array}
\right.
\end{equation}
Then $u_\varepsilon \in C^2[0,L]$, the function $a_\varepsilon(x)
\equiv \frac{-u_\varepsilon ''(x)}{u_\varepsilon (x)}, \ \forall \
x \in [0,L], \ x \neq \frac{(2k-1)L}{2(n+1)},  \ 1 \leq k \leq
n+1,$ belongs to $\landan$ and
\begin{equation}\label{0307t2}\liminf_{\varepsilon \rightarrow
0^+} \ \Vert a_\varepsilon - \landaene \Vert_{L^1 (0,L)} =
\frac{n\pi}{L} 2(n+1) \cot \frac{n\pi}{2(n+1)}.
\end{equation}
\end{lemma}
\begin{proof}
We claim that for each $0 \leq i \leq 2n+1,$ function
$a_\varepsilon$ satisfies
\begin{equation}\label{0307t1}
\landaene \prec a_\varepsilon, \ \mbox{in the interval } \ \left (
\frac{iL}{2(n+1)}, \frac{(i+1)L}{2(n+1)} \right )
\end{equation}
and
\begin{equation}\label{0307t3}
\liminf_{\varepsilon \rightarrow 0^+} \ \Vert a_\varepsilon -
\landaene \Vert_{L^1 (\frac{iL}{2(n+1)}, \frac{(i+1)L}{2(n+1)})} =
\frac{n\pi}{L} \cot \frac{n\pi}{2(n+1)}
\end{equation}
It is trivial that from (\ref{0307t1}) and (\ref{0307t3}) we
deduce (\ref{0307t2}). Moreover, taking into account the
definition of the function $u_\varepsilon,$ it is clear that it is
sufficient to prove the claim in the case $i =0.$ Now, if $x \in
(0,\frac{L}{2(n+1)})$ we can distinguish two cases:
\begin{enumerate}
\item $x \in (\varepsilon, \frac{L}{2(n+1)}).$ Then $a_\varepsilon
(x) = \frac{-u_\varepsilon ''(x)}{u_\varepsilon (x)} \equiv
\landaene.$ \item $x \in (0,\varepsilon).$ Then
\end{enumerate}
$$
a_\varepsilon (x) - \landaene = \frac{-2
\frac{x-\varepsilon}{\varepsilon^2} \frac{n\pi}{L} \cos
\frac{n\pi}{2(n+1)} -  \frac{(x-\varepsilon)^3}{3\varepsilon^2}
\frac{n^3\pi^3}{L^3} \cos \frac{n\pi}{2(n+1)}}{-\sin
(\frac{n\pi}{L}(x-\frac{L}{2(n+1)})) +
\frac{(x-\varepsilon)^3}{3\varepsilon^2} \frac{n\pi}{L} \cos
\frac{n\pi}{2(n+1)}} > 0
$$
Therefore $a_\varepsilon \in \Lambda_n.$ Moreover, if $\varepsilon
\rightarrow 0^+,$ then
$$
\frac{- \frac{(x-\varepsilon)^3}{3\varepsilon^2}
\frac{n^3\pi^3}{L^3} \cos \frac{n\pi}{2(n+1)}}{-\sin
(\frac{n\pi}{L}(x-\frac{L}{2(n+1)})) +
\frac{(x-\varepsilon)^3}{3\varepsilon^2} \frac{n\pi}{L} \cos
\frac{n\pi}{2(n+1)}} \rightarrow 0,
$$
uniformly if $x \in (0,\varepsilon).$

Finally, since
$$
\lim_{\varepsilon \rightarrow 0^+} \int_0^{\varepsilon} \left [
\frac{-2 \frac{x-\varepsilon}{\varepsilon^2} \frac{n\pi}{L} \cos
\frac{n\pi}{2(n+1)}} {-\sin (\frac{n\pi}{L}(x-\frac{L}{2(n+1)})) +
\frac{(x-\varepsilon)^3}{3\varepsilon^2} \frac{n\pi}{L} \cos
\frac{n\pi}{2(n+1)}} - \frac{-2
\frac{x-\varepsilon}{\varepsilon^2} \frac{n\pi}{L} \cos
\frac{n\pi}{2(n+1)}} {-\sin (\frac{n\pi}{L}(x-\frac{L}{2(n+1)})) }
\right ] = 0
$$
and
$$-\sin (\frac{n\pi}{L}(x-\frac{L}{2(n+1)})) \rightarrow \sin
\frac{n\pi}{2(n+1)},$$ uniformly  in $x \in (0,\varepsilon)$ when
$\varepsilon \rightarrow 0+,$ we deduce
$$
\liminf_{\varepsilon \rightarrow 0^+} \ \Vert a_{\varepsilon} -
\landaene \Vert_{L^1 (0,\frac{L}{2(n+1)})} = \liminf_{\varepsilon
\rightarrow 0^+} \frac{n\pi}{L} \cot \frac{n\pi}{2(n+1)}
\frac{2}{\varepsilon^2}\int_0^{\varepsilon} (\varepsilon - x) =
\frac{n\pi}{L} \cot \frac{n\pi}{2(n+1)}
$$
which is (\ref{0307t3}) for the case $i=0$.
\end{proof}

\begin{lemma}\label{l6}
$\beta_{1,n}$ is not attained.
\end{lemma}

\begin{proof}
Let $a \in \landan$ be such that $\Vert a - \landaene \Vert_{L^1
(0,L)} = \beta_{1,n}.$ Let $u$ be any nontrivial solution of
(\ref{n1}) associated to the function $a$. As previously, we
denote the zeros of $u'$ by $0 = x_0 < x_2 < \ldots < x_{2m} = L$
and the zeros of $u$ by $ x_1 < x_3 < \ldots < x_{2m-1}. $ By
using Lemma \ref{l3}, Lemma \ref{l4} and Lemma \ref{l7},  we have
\begin{equation}\label{1209071}
\begin{array}{c}
\beta_{1,n} = \Vert a - \landaene \Vert_{L^1 (0,L)} =
\displaystyle \sum_{i=0}^{2m-1}  \Vert a - \landaene \Vert_{L^1
(x_i,x_{i+1})} \geq \\ \\
\sum_{i=0}^{2m-1} J_i (u) \geq \frac{n\pi}{L}\displaystyle
\sum_{i=0}^{2m-1} \cot \frac{n\pi (x_{i+1}-x_i)}{L} \geq \\ \\
\frac{n\pi}{L} 2m \cot \frac{n\pi}{2m} \geq \frac{n\pi}{L} 2(n+1)
\cot \frac{n\pi}{2(n+1)} = \beta_{1,n}
\end{array}
\end{equation}
where $J_i (u)$ is given either by $$ J_i(u) =
\frac{\int_{x_i}^{x_{i+1}} u'^2 - \landaene \int_{x_i}^{x_{i+1}}
u^2 }{u^2(x_{i+1})}, \ \ \mbox{if} \ u(x_i) = 0 $$ or by
$$
J_i(u)= \frac{\int_{x_i}^{x_{i+1}} u'^2 - \landaene
\int_{x_i}^{x_{i+1}} u^2 }{u^2(x_i)}, \ \ \mbox{if} \ u(x_{i+1}) =
0. $$ Consequently, all inequalities in (\ref{1209071}) transform
into equalities. In particular we obtain from Lemma \ref{l4} and
the property (P) shown in Lemma \ref{l5} that
$$
m = n+1, \ x_{i+1}-x_i = \frac{L}{2(n+1)}, \ 0 \leq i \leq 2n+1.
$$
Also, it follows
$$ J_i (u) = \frac{n\pi}{L} \cot \frac{n\pi}{ L}\frac{L}{ 2 (n+1)}, \
0 \leq i \leq 2n+1.
$$
From Lemma \ref{l2} we deduce that, up to some nonzero constants,
function $u$ fulfils in each interval $[x_i,x_{i+1}],$
$$
u(x) = \frac{\sin \frac{n\pi}{L} (x-x_i)}{\sin \frac{n\pi}{L}
(x_{i+1}-x_i)}, \ \mbox{if} \ i \ \mbox{is odd},
$$
$$
\mbox{and}\ \ u(x) = \frac{\sin \frac{n\pi}{L} (x-x_{i+1})}{\sin
\frac{n\pi}{L} (x_{i}-x_{i+1})}, \ \mbox{if} \ \ i \ \mbox{is
even}.
$$
In particular, in the interval $[0,\frac{L}{2(n+1)}],$ $u$ must be
the function
$$
u(x) = \frac{\sin \frac{n\pi}{L}(x-\frac{L}{2(n+1)})}{\sin
\frac{n\pi}{L}(-\frac{L}{2(n+1)})}
$$
which does not satisfy the condition $u'(0) = 0.$ The conclusion
is that $\beta_{1,n}$ is not attained.

\end{proof}

Finally, as a trivial consequence of Lemma \ref{l5}, Lemma
\ref{l7} and Lemma \ref{l6} we have the conclusion of Theorem
\ref{t1}.

\end{proof}

\begin{remark}
Let us observe that if we consider $\beta_{1,n}$ as a function  of
$n \in (0,+\infty),$ then $\lim_{n \rightarrow 0^+} \ \beta_{1,n}
= \frac{4}{L},$ the constant of the classical $L^1$ Lyapunov
inequality at the first eigenvalue (\cite{ha}).
\end{remark}
\begin{remark}\label{1212073}
The case where $L = 1$ and function $a$ satisfies the condition $
A \leq a(x) \leq B, \ \mbox{a.e. in} \ (0,L) $ where $\lambda_k <
A < \lambda_{k+1} \leq B$ for some $k \in \natu \cup \{ 0 \}, $
has been considered in \cite{yong}, where the authors use Optimal
Control theory methods. In this paper, the authors define the set
$\Lambda_{A,B}$ as the set of functions $a \in L^{1}(0,L)$ such
that $A \leq a(x) \leq B, \ \mbox{a.e. in} \ (0,L)$ and (\ref{n1})
has nontrivial solutions. Then, by using the Pontryagin's maximum
principle they prove that the number
$$ \beta_{A,B} \equiv \inf_{a \in
\Lambda_{A,B} } \ \Vert a \Vert_{L^1 (0,L)} $$ is attained. In
addition, they calculate $\lim_{B \rightarrow +\infty} \
\beta_{A,B}.$
\end{remark}
\begin{remark} In our opinion, the inequality
$\int_0^1 b(t) \ dt \leq 2 \sqrt{A}\cot \frac{\sqrt{A}}{2}$ in
\cite{yong}, Theorem 3, must be substituted by $\int_0^1 b(t) \ dt
\leq A + 2(k+1) \sqrt{A}\cot \frac{\sqrt{A}}{2(k+1)}.$ This may be
easily derived from our method modifying the definition of the set
$\Lambda_n$ (given in (\ref{2806073})) in a trivial way.
\end{remark}
\begin{remark}
If $A \rightarrow \lambda{_k}^{+},$ it does not seem possible to
deduce from \cite{yong} that the constant $\beta_{1,k}$ (defined
in (\ref{1009071})) is not attained. In fact, to the best of our
knowledge, this result is new. Moreover, our method, which
combines a detailed analysis about the number and distribution of
zeros of nontrivial solutions of (\ref{n1}) and their first
derivatives, together with the use of suitable minimization
problems, will be very useful to combine Lyapunov inequalities and
disfocality. This will be seen in the next section.
\end{remark}
\begin{remark}
We can use our methods to do an analogous study for other boundary
conditions. In particular with the help of Lemma \ref{280607l1}
and Lemma \ref{l2} we can consider the mixed linear problem
\begin{equation}\label{1212071}
u''(x) + a(x)u(x) = 0, \ x \in (0,L), \ u'(0) = u(L) = 0
\end{equation}
where
$$ a \in \Gamma_n = \{ a \in L^1 (0,L): \mu_n \prec a \
\mbox{and (\ref{1212071}) has nontrivial solutions} \}$$ Here
$\mu_n$ is the $n$-th eigenvalue of the eigenvalue problem
\begin{equation}
 u''(x) + \mu u(x) = 0, \ x \in (0,L), \ u'(0) = u(L) =
0 \end{equation}
\end{remark}
The case where $L = 1$ and function $a$ satisfies the condition $
A \leq a(x) \leq B, \ \mbox{a.e. in} \ (0,L) $ where $\mu_k < A <
\mu_{k+1} \leq B$ has been considered in \cite{yuhuayong}. As in
\cite{yong}, the authors use Optimal Control theory methods. See
also \cite{huai} for Dirichlet boundary conditions.

\section{Lyapunov inequalities and disfocality}

The $L^\infty$ Lyapunov inequality is trivial from Dolph's result
(\cite{dolp}). In fact, by using Dolph's result, the constant
\begin{equation}\label{041207} \beta_{\infty,n} \equiv \inf_{a \in \landan } \
\Vert a \Vert_{L^\infty (0,L)}
\end{equation}
must be greater than or equal to $\lambda_{n+1}.$ Since the
constant function $\lambda_{n+1}$ is an element of $\landan,$ we
deduce
\begin{equation}
\beta_{\infty,n} = \lambda_{n+1}.
\end{equation}
Moreover $\beta_{\infty,n}$ is attained in a unique element
$a_{\infty} \in \landan$ given by the constant function $a_\infty
\equiv \lambda_{n+1}$.

\medskip

On the other hand, under the restriction
\begin{equation}\label{2211073} a \in L^1 (0,L), \ \landaene \prec
a, \end{equation} the relation between Neumann boundary conditions
and disfocality arises in a natural way. In fact, if $u \in \sobol
$ is any nontrivial solution of (\ref{n1}) and the zeros of $u$
are denoted by $ x_1 < x_3 < \ldots < x_{2m-1},$ and the zeros of
$u'$ are denoted by $0 = x_0 < x_2 < \ldots < x_{2m} = L,$ then
for each given $i, \ 0 \leq i \leq 2m-1,$ function $u$ satisfies
{\small
\begin{equation}\label{2211071} u''(x) + a(x)u(x) = 0, x \in
(x_i,x_{i+1}), u(x_i) = 0, u'(x_{i+1}) = 0,\ \ \mbox{if} \ i \
\mbox{is odd }
\end{equation} and
\begin{equation}\label{2211072} u''(x) + a(x)u(x) = 0, \ x \in (x_i,x_{i+1}),
\ \ u'(x_i) = 0, \ u(x_{i+1}) = 0, \ \mbox{if} \ i \ \mbox{is even
}.
\end{equation}
} In consequence, each one of the problems (\ref{2211071}) and
(\ref{2211072}) with $0 \leq i \leq 2m-1,$ have nontrivial
solution. This simple observation can be used to deduce the
following conclusion: if $a$ is any function satisfying
(\ref{2211073}) such that for any $m \geq n+1$ and any
distribution of numbers $0 = x_0 < x_1 < x_2< \ldots < x_{2m-1} <
x_{2m} = L,$ either some problem of the type (\ref{2211071}) or
some problem of the type (\ref{2211072}) has only the trivial
solution, then problem (\ref{n1}) has only the trivial solution.
Lastly, it has been established in \cite{cavi} (Theorem 2.1 for
the case $p = \infty$) that if $b \in L^\infty (c,d)$ satisfies
\begin{equation}\label{2211074}
\Vert b \Vert_{L^{\infty}(c,d)} \leq \frac{\pi^2}{4(d-c)^2} \
\mbox{and} \ b \neq \frac{\pi^2}{4(d-c)^2} \ \mbox{in} \ (c,d)
\end{equation}
then the unique solution of the boundary value problems
\begin{equation}\label{m22bis}
u''(x) + b(x)u(x) = 0, \ x \in (c,d), \ u'(c) = u(d) = 0
\end{equation} and
\begin{equation}\label{m22bisbis}
u''(x) + b(x)u(x) = 0, \ x \in (c,d), \ u(c) = u'(d) = 0
\end{equation}
is the trivial one.

We may use previous reasonings to obtain the following result
\begin{theorem}\label{t2211}
If function $a$ fulfils
\begin{equation}\label{2211076}
\begin{array}{c}
a \in L^\infty (0,L), \ \landaene \prec a \
\mbox{and} \ \exists \ \ 0 = y_0 < y_1 < \ldots < y_{2n+1} < y_{2n+2} = L: \\ \\
\max\ _{0 \leq i \leq 2n+1}  \{ (y_{i+1}-y_i)^2 \Vert a
\Vert_{L^\infty (y_i,y_{i+1})} \} \leq \pi^2/4
\end{array}
\end{equation}
 and, in addition, $a$ is not the
constant $\pi^2 /4(y_{i+1}-y_i)^2,$ at least in one of the
intervals $[y_i,y_{i+1}], \ 0 \leq i \leq 2n+1,$

then the boundary value problem (\ref{n1}) has only the trivial
solution.
\end{theorem}

\begin{proof}
To prove this Theorem, take into account that if $m \geq n+1$ and
$0 = x_0 < x_1 < x_2< \ldots < x_{2m-1} < x_{2m} = L,$ is any
arbitrary distribution of numbers, then or
\begin{equation}\label{2211078}
[x_j,x_{j+1}] \subset [y_i,y_{i+1}], \ \ \mbox{strictly,}  \
\end{equation}
for some $0 \leq i \leq 2n+1, \ 0 \leq j \leq 2m-1$ or
\begin{equation}\label{22110710}
m = n+1 \ \mbox{and} \ x_i = y_i, \ \forall \ 0 \leq i \leq \
2n+2.
\end{equation}
If (\ref{2211078}) is satisfied, then
\begin{equation}\label{2211079}
\Vert a \Vert_{L^\infty (x_j,x_{j+1})} < \Vert a \Vert_{L^\infty
(y_i,y_{i+1})} \leq \frac{\pi^2}{4(y_{i+1}-y_i)^2} <
\frac{\pi^2}{4(x_{j+1}-x_j)^2}
\end{equation}
and consequently we deduce from (\ref{2211071}), (\ref{2211072})
and (\ref{2211074}) that (\ref{n1}) has only the trivial solution.

If (\ref{22110710}) is satisfied, we deduce from the hypotheses of
the Theorem, that  $a$ is not the constant $\pi^2
/4(x_{i+1}-x_i)^2,$ at least in one of the intervals
$[x_i,x_{i+1}], \ 0 \leq i \leq 2n+1.$ Therefore, again
(\ref{2211071}), (\ref{2211072}) and (\ref{2211074}) imply that
(\ref{n1}) has only the trivial solution. In any case, we have the
desired conclusion.
\end{proof}

\begin{remark}
If in previous Theorem we choose $y_i = \frac{iL}{2(n+1)},\ 0 \leq
i \leq 2n+2,$ then we have the so called non-uniform non-resonance
conditions at higher eigenvalues (\cite{dolp}, \cite{mawa3}) but
if for instance, $y_{j+1}-y_{j} < \frac{L}{2(n+1)},$ for some $j,
\ 0 \leq j \leq 2n+1,$ function $a$ can satisfies $\Vert a
\Vert_{L^\infty (y_j,y_{j+1})} = \frac{\pi^2}{4(y_{j+1}-y_j)^2}$
(which is a quantity greater than $\lambda_{n+1} =
\frac{(n+1)^2\pi^2}{L^2}$) as long as $a$ satisfies
(\ref{2211076}) for each $i \neq j.$
\end{remark}
\begin{remark}\label{r221107}
The hypothesis of the previous Theorem is optimal in the sense
that if $a$ is the constant $\pi^2 /4(y_{i+1}-y_i)^2$ in each one
of the intervals $(y_i,y_{i+1}), \ 0 \leq i \leq 2n+1,$ then
(\ref{n1}) has nontrivial solutions. In fact, if this is the case,
it is easily checked that there exist appropriate constants $k_i,
\ 0 \leq i \leq 2n+1,$ such that the function
$$
u(x) = \left \{ \begin{array}{l} k_i \cos
\frac{\pi(x-y_i)}{2(y_{i+1}-y_i)}, \ x \in [y_i,y_{i+1}], \  i\
\mbox{even}, \\ \\
 k_i \cos
\frac{\pi(y_{i+1}-x)}{2(y_{i+1}-y_i)}, \ x \in [y_i,y_{i+1}], \ i\
\mbox{odd},
\end{array}
\right.
$$
\end{remark}
is a nontrivial solution of (\ref{n1}).

\medskip

Now we comment some relations between the Lyapunov constant
$\beta_{1,n},$ given in Theorem \ref{t1} and disfocality. To this
respect, it is clear from the definition of $\beta_{1,n},$ that if
a function $a$  satisfies
\begin{equation}\label{2311071}
a \in L^1(0,L), \ \lambda_n \prec a, \ \Vert a - \lambda_n \Vert_1
< \beta_{1,n}
\end{equation}
then the unique solution of (\ref{n1}) is the trivial one. In the
next Theorem we prove that, with the use of disfocality, we can
obtain a more general condition.

\begin{theorem}\label{t111207}
$ $
\begin{enumerate}
\item If function $a \in L^1(0,L), \ \lambda_n \prec a, $
satisfies:
\begin{equation}\label{2311074}
\begin{array}{c}
\exists \  \ 0 = y_0< y_1<\ldots < y_{2n+1} <
y_{2n+2} = L:  \\ \\
y_{i+1}-y_i<\frac{L}{2n}\, ;\ \Vert a - \lambda_n \Vert
_{L^1(y_{i},y_{i+1})} < \frac{n\pi}{L} \cot
\frac{n\pi(y_{i+1}-y_i)}{L}, \ \forall \ 0 \leq i \leq 2n+1,
\end{array}
\end{equation}
then the unique solution of (\ref{n1}) is the trivial one. \item
(\ref{2311071}) implies (\ref{2311074}). \item If $0 = y_0 < y_1 <
\ldots < y_{2n+1} < y_{2n+2} = L,$ is any distribution of numbers
such that $y_{k+1}-y_k<\frac{L}{2n}, \ \forall \ 0 \leq k \leq
2n+1$ and $y_{i+1}-y_i \neq y_{j+1}-y_j,$ for some $0 \leq i,j
\leq 2n+1,$ then there exists $a \in L^{1}(0,L), \ \lambda_n \prec
a,$ satisfying (\ref{2311074}) but not satisfying (\ref{2311071}).
\end{enumerate}
\end{theorem}

\begin{proof}
If $a$ satisfies (\ref{2311074}), then the unique solution of
(\ref{n1}) is the trivial one. In fact, if this is not true, let
$u$ be a nontrivial solution of (\ref{n1}) and let us denote the
zeros of $u$ by $ x_1 < x_3 < \ldots < x_{2m-1}$ and the zeros of
$u'$ by $0 = x_0 < x_2 < \ldots < x_{2m} = L.$ Since $m \geq n+1,$
then
\begin{equation}\label{2311075}
[x_j,x_{j+1}] \subset [y_i,y_{i+1}]
\end{equation}
for some $0 \leq i \leq 2n+1, \ 0 \leq j \leq 2m-1.$ Consequently,
$$
\displaystyle \frac{\Vert a - \lambda_n
\Vert_{L^1(x_j,x_{j+1})}}{\cot \frac{n\pi(x_{j+1}-x_{j})}{L}} \leq
\frac{\Vert a - \lambda_n \Vert_{L^1(y_i,y_{i+1})}}{\cot
\frac{n\pi(y_{i+1}-y_{i})}{L}} < \frac{n\pi}{L}.
$$
From here we deduce
$$
\Vert a - \lambda_n \Vert_{L^1(x_j,x_{j+1})} < \frac{n\pi}{L}\cot
\frac{n\pi(x_{j+1}-x_{j})}{L}
$$
which is a contradiction with Lemma \ref{280607l1} and Lemma
\ref{l2}.

Next we prove that (\ref{2311071}) implies (\ref{2311074}). We can
certainly assume that $\inf a>\lambda_n$, for if not, we replace
$a$ by $a +\delta$ (for small $\delta>0$) and the new function
$a+\delta$ satisfies (\ref{2311071}). Note that if condition
(\ref{2311074}) is satisfied for $a+\delta$ then also is satisfied
for the function $a$.

Now choose $\varepsilon> 0$ sufficiently small. Since the function
$$
 \displaystyle \frac{\Vert a - \lambda_n
\Vert _{L^1(0,y)}}{\cot \frac{n\pi(y-0)}{L}}
$$
is strictly increasing with respect to $y \in (0,\frac{L}{2n})$
and
$$
\lim_{y \rightarrow 0^+} \ \displaystyle \frac{\Vert a - \lambda_n
\Vert _{L^1(0,y)}}{\cot \frac{n\pi(y-0)}{L}} = 0, \ \ \lim_{y
\rightarrow \frac{L}{2n}^-} \ \displaystyle \frac{\Vert a -
\lambda_n \Vert _{L^1(0,y)}}{\cot \frac{n\pi(y-0)}{L}} = +\infty
$$
there is an unique $y_1,$  $0 = y_0 < y_{1} < \frac{L}{2n}$ such
that
\begin{equation}\label{1112072}
 \displaystyle \frac{\Vert a - \lambda_n
\Vert _{L^1(0,y_1)}}{\cot \frac{n\pi(y_1 -0)}{L}} = \frac{n\pi}{L}
- \varepsilon.
\end{equation}
With the help of a similar reasoning, it is possible to prove the
existence of points $0 = y_0 < y_1 < \ldots < y_{2n+1},$ such that
\begin{equation}\label{1112073}
 \displaystyle \frac{\Vert a - \lambda_n
\Vert _{L^1(y_i,y_{i+1})}}{\cot \frac{n\pi(y_{i+1} -y_i)}{L}} =
\frac{n\pi}{L} - \varepsilon, \ \ y_{i+1} - y_i < \frac{L}{2n}, \
\  0 \leq i \leq 2n.
\end{equation}
(If it is necessary, we can define $a(x) = \lambda_n, \ \forall x
> L$).

Since $y_{i+1}-y_i<\frac{L}{2n}, \ 0 \leq i \leq 2n-1$, then
$y_{2n}<L$.

If $y_{2n+1}\geq L$, then we replace the number $y_{2n+1}$ with
$y_{2n+1}=L-\mu$ (for small $\mu>0$). Finally, choosing $y_{2n+2}
= L$, we obtain (\ref{2311074}).

If $y_{2n+1}<L$, take $y_{2n+2} = L$. We claim that
\begin{equation}\label{1112074}
y_{2n+2}-y_{2n+1}<\frac{L}{2n} \ \ \mbox{and} \ \  \displaystyle
\frac{\Vert a - \lambda_n \Vert _{L^1(y_{2n+1},y_{2n+2})}}{\cot
\frac{n\pi(y_{2n+2} -y_{2n+1})}{L}} <  \frac{n\pi}{L} -
\varepsilon.
\end{equation}
In fact, if $y_{2n+2}-y_{2n+1}\geq \frac{L}{2n},$ then
$y_{2n+1}\leq \frac{L(2n-1)}{2n}$. Then, from (\ref{1112073}),
Lemma \ref{l4} (with $r=2n+1$, $S=\frac{n\pi}{L}(y_{2n+1})$ and
$z_i=\frac{n\pi}{L}(y_{i+1}-y_i)$) and using the monotonicity of
$\cot$ in $(0,\pi/2)$ we obtain
$$
\begin{array}{c}
\displaystyle \frac{n\pi}{L}2(n+1)\cot\,
\frac{n\pi}{2(n+1)}=\beta_{1,n}>\sum_{i=0}^{2n}\Vert a - \lambda_n
\Vert_{L^1 (y_i,y_{i+1})} = \\ \\ \displaystyle \left(
\frac{n\pi}{L}-\varepsilon \right)\sum_{i=0}^{2n}\cot
\frac{n\pi}{L}(y_{i+1}-y_i)\geq \left( \frac{n\pi}{L}-\varepsilon
\right) (2n+1)\cot \frac{n\pi}{L(2n+1)}y_{2n+1}\geq
\\ \\ \displaystyle \left( \frac{n\pi}{L}-\varepsilon
\right) (2n+1)\cot \frac{\pi(2n-1)}{2(2n+1)}
\end{array}
$$
If $\varepsilon \rightarrow 0^+,$ we conclude
\begin{equation}\label{2012071}
\beta_{1,n} \geq \displaystyle \frac{n\pi}{L} (2n+1)\cot
\frac{\pi(2n-1)}{2(2n+1)}
\end{equation}

Now, by using that the function $\displaystyle x\mapsto \frac{2\pi
\cot x}{\pi-2x}$ is strictly decreasing in $(0, \pi/2)$ and that
$\displaystyle \frac{\pi(2n-1)}{2(2n+1)} <\frac{n\pi}{2(n+1)}$, we
obtain
$$
\beta_{1,n} \geq  \displaystyle \frac{n\pi}{L} (2n+1)\cot
\frac{\pi(2n-1)}{2(2n+1)} > \displaystyle
\frac{n\pi}{L}2(n+1)\cot\, \frac{n\pi}{2(n+1)} = \beta_{1,n}
$$
\noindent which is a contradiction.

\noindent It remains to prove the second part of the claim
(\ref{1112074}). In fact, if this second part is not true, then
from (\ref{1112073}) and Lemma \ref{l4} (with $r=2n+2$, $S=n\pi$
and $z_i=\frac{n\pi}{L}(y_{i+1}-y_i)$) we have
$$
\begin{array}{c}
\Vert a - \lambda_n \Vert_{L^1(0,L)} = \displaystyle
\sum_{i=0}^{2n+1} \Vert
a - \lambda_n \Vert_{L^1 (y_i,y_{i+1})} \geq  \\ \\
(\frac{n\pi}{L} - \varepsilon ) \displaystyle \sum_{i=0}^{2n+1}
\cot \frac{n\pi(y_{i+1} -y_i)}{L} \geq (\frac{n\pi}{L} -
\varepsilon ) \displaystyle \frac{\beta_{1,n}}{n\pi/L},
\end{array}
$$
for $\varepsilon > 0$ sufficiently small. This is a contradiction
with (\ref{2311071}).

Finally, to prove part (3) of the theorem, let us take numbers $0
= y_0 < y_1 < \ldots < y_{2n+1} < y_{2n+2} = L,$ such that
$y_{k+1}-y_k<\frac{L}{2n}, \ \forall \ 0 \leq k \leq 2n+1$ and
$y_{i+1}-y_i \neq y_{j+1}-y_j,$ for some $0 \leq i,j \leq 2n+1.$
Then from Lemma \ref{l4} we obtain
$$
\sum_{i=0}^{2n+1} \frac{n\pi}{L} \cot\frac{n\pi(y_{i+1}-y_i)}{L} >
\frac{2\pi n(n+1)}{L}\cot\frac{n\pi}{2(n+1)} = \beta_{1,n}
$$
Now, choose a function $a \in L^{1}(0,L), \ \lambda_n \prec a,$
satisfying
$$
\Vert a - \lambda_n \Vert _{L^1(y_{i},y_{i+1})} = \frac{n\pi}{L}
\cot \frac{n\pi(y_{i+1}-y_i)}{L} - \varepsilon, \ \forall \ 0 \leq
i \leq 2n+1
$$
It is trivial that if $\varepsilon$ is sufficiently small, then
function $a$ satisfies (\ref{2311074}) whereas
$$
\Vert a - \lambda_n \Vert _{L^{1}(0,L)} = \sum_{i=0}^{2n+1} \Vert
a - \lambda_n \Vert _{L^{1}(y_i,y_{i+1})} > \beta_{1,n}.
$$
\end{proof}

\medskip

{\bf Final remark on nonlinear problems.}

We finish this paper by showing how to use previous reasonings to
obtain new theorems on the existence and uniqueness of solutions
of nonlinear b.v.p.
\begin{equation}\label{m38}
u''(x) + f(x,u(x)) = 0, \ x \in (0,L), \ \ u'(0) = u'(L) = 0.
\end{equation}

For example, we have the following theorem related to Theorem 2.1
in \cite{mawa3}. This last Theorem allows to consider more general
boundary value problems, but for ordinary problems with Neumann
boundary conditions our hypotheses allow a more general behavior
on the derivative $f_u (x,u).$ We omit the details of the proof
(see \cite{camovi} and \cite{camovi2} for similar results at the
two first eigenvalues).
\begin{theorem}\label{t2bis}
Let us consider (\ref{m38}) where the following requirements are
supposed:
\begin{enumerate}
\item $f$ and $f_{u}$ are Caratheodory functions on $[0,L]\times
\mathbb R$ and $f(\cdot,0) \in L^1 (0,L).$ %
\item There exist functions $\alpha, \beta \in L^{\infty} (0,L),$
satisfying
$$\lambda_n \leq \alpha (x) \leq f_{u}(x,u) \leq \beta (x)$$ on $[0,L]\times \mathbb
R.$ Furthermore, $\alpha$ differs from $\lambda_n $ on a set of
positive measure and $\beta$ satisfies either hypothesis
(\ref{2211076}) of Theorem \ref{t2211} or hypothesis
(\ref{2311074}) of Theorem \ref{t111207}.
\end{enumerate}
Then, problem (\ref{m38}) has a unique solution.
\end{theorem}

\end{document}